\documentclass[fleqn]{amsart}
\usepackage{amssymb, graphics, amsmath, epsfig, changebar, hyperref}
\usepackage{enumitem}
\usepackage{mathtools} 
\usepackage[percent]{overpic}
\usepackage[caption=false]{subfig} 

\theoremstyle{plain}
\newtheorem{theorem}{Theorem}
\newtheorem{corollary}{Corollary}

\newtheorem{lemma}{Lemma}
\newtheorem{proposition}{Proposition}

\newtheorem{question}{Question}

\theoremstyle{definition}

\theoremstyle{example}

\theoremstyle{remark}

\numberwithin{equation}{section}
\newcommand{\N}{\mathbb{N}}

\newcommand{\Z}{\mathbb{Z}}

\newcommand{\C}{\mathbb{C}}

\newcommand{\Sum}{\displaystyle\sum}

\newcommand{\Om}{\Omega}

\newcommand{\Phit}{\widetilde{\Phi}}
\newcommand{\No}{\N\cup\{0\}}

\graphicspath{{./images/}}

\newcommand{\overimage}[3]{\begin{overpic}{#1}\put(#2){#3}\end{overpic}}
\newcommand{\oversimage}[5]{\begin{overpic}{#1}\put(#2){#3}\put(#4){#5}\end{overpic}}
\newcommand{\overssimage}[7]{\begin{overpic}{#1}\put(#2){#3}\put(#4){#5}\put(#6){#7}\end{overpic}}

\begin{document}                   

\title{Infinitely many roots of unity are zeros of some Jones polynomials}
\author{Maciej Mroczkowski}
\address{Institute of Mathematics\\
Faculty of Mathematics, Physics and Informatics\\
University of Gdansk, 80-308 Gdansk, Poland\\
e-mail: maciej.mroczkowski@ug.edu.pl}

\begin{abstract}
Let $N=2n^2-1$ or $N=n^2+n-1$, for any $n\ge 2$. Let $M=\frac{N-1}{2}$.
We construct families of prime knots with Jones polynomials $(-1)^M\sum_{k=-M}^{M} (-1)^kt^k$.  Such polynomials have Mahler measure equal to $1$. If $N$ is prime, these are cyclotomic polynomials $\Phi_{2N}(t)$, up to some shift in the powers of $t$. Otherwise, they are products of such polynomials, including $\Phi_{2N}(t)$. In particular, all roots of unity $\zeta_{2N}$ occur as roots of Jones polynomials. We also show that some roots of unity cannot be zeros of Jones polynomials.
\end{abstract}
\maketitle

\let\thefootnote\relax\footnotetext{Mathematics Subject Classification 2020: 57K10, 57K14} 

\section{Introduction}
We study knots $K$ with Jones polynomials that have Mahler measure equal to $1$, or $M(V_{K}(t))=1$.
Such knots were considered in~\cite{CK,CK2}. A Laurent polynomial $P$, with $M(P)=1$, has the form: $t^a$ times a product of cyclotomic polynomials $\Phi_n$, $a\in\Z$. Following~\cite{CK2}, we call such $P$ {\it cyclotomic}.

The motivation for studying such knots comes from some observed connections between the Mahler measure of the Jones polynomial of a knot and its hyperbolic volume~\cite{CK,CK2}. The Mahler measure of Jones polynomials has also been studied in~\cite{T1, T2}. In a more general context, not much is known about the question: what polynomials are Jones polynomials? This is in contrast to the Alexander polynomial: there are simple conditions on a polynomial which are sufficient and necessary for it to be the Alexander polynomial of a knot. Studying the locus of the zeros of Jones polynomials is part of the general question.
It is shown in \cite{JZDT} that this locus is dense in $\C$. Our result implies that the locus intersected with the unit circle is dense in the unit circle. 

Near the end of~\cite{CK}, after listing all knots up to $16$ crossings with cyclotomic Jones polynomials (there are only $17$ such knots), the following problem is posed:
``An interesting open question is how to construct more knots with $M(V_K(t))=1$.''
In this paper, we construct four inifinite families of knots with cyclotomic Jones polynomials of a particularly simple form.
They include $4_1$ and $9_{42}$ with $V_{4_1}(t)=t^{-2}-t^{-1}+1-t+t^2$ and $V_{9_{42}}(t)=t^{-3}-t^{-2}+t^{-1}-1+t-t^2+t^3$. The Jones polynomials of the other knots extend these two examples. Their coefficients are finite sequences of alternating $1$'s and $-1$'s, starting and ending with $1$. In particular, there is no bound on the span of such Jones polynomials.

In fact, polynomials with alternating $1$ and $-1$'s as coefficients, starting and ending with $1$ are obtained as follows. Let $m$ be odd. It is well known that, for any $a\in N$, $t^a-1=\prod_{k|a}\Phi_k(t)$.
Hence, 
\[t^{2m}-1=\prod_{k|2m}\Phi_k(t)=\prod_{k|m}\Phi_k(t)\prod_{k|m}\Phi_{2k}(t)=(t^m-1)\prod_{k|m}\Phi_{2k}(t)\]
It follows, that $\prod_{k|m}\Phi_{2k}(t)=t^m+1=(t+1)(t^{m-1}-t^{m-2}+t^{m-3}+\ldots-t+1)$.
Since $\Phi_2(t)=t+1$, we get:
\[\prod_{k|m,k>1}\Phi_{2k}(t)=t^{m-1}-t^{m-2}+t^{m-3}+\ldots-t+1\]
For $m$ odd, we introduce the following notation:
\[\Phit_{2m}(t):=t^{-\frac{m-1}{2}}\prod_{k|m,k>1}\Phi_{2k}(t)=t^{-\frac{m-1}{2}}-t^{-\frac{m-1}{2}+1}+\ldots-
t^{\frac{m-1}{2}-1}+t^{\frac{m-1}{2}}\]
We allow $m=1$ with $\Phit_2(t)=1$.
We will construct knots with cyclotomic Jones polynomials equal to $\Phit_{2m}$ for infinitely many odd $m$.

Notice, that $\frac{d}{dt}[t^nP(t)]_{t=1}=P'(1)+n$ for any $n\in\Z$,
if $P$ is a Laurent polynomial satisfying $P(1)=1$. If $P$ is a Jones polynomial of a knot, it satisfies $P(1)=1$ and $P'(1)=0$ (see \cite{J}),
hence $t^nP$ cannot be a Jones polynomial for $n\neq 0$. 
We say that a Laurent polynomial $P$ is {\it palindromic}, if $P(t^{-1})=t^nP(t)$, for some $n\in \Z$; in particular, if $n=0$,
we say that it is {\it symmetric}. One checks that $P'(1)=0$, if $P$ is symmetric. Hence, a palindromic Jones polynomial of a knot must be symmetric.

For $n\ge 3$, cyclotomic polynomials $\Phi_n$ are palindromic of even degree, but not symmetric (since they are not Laurent polynomials).  In order to make them symmetric, we multpily $\Phi_n$ with $t^{-\frac{\varphi(n)}{2}}$ (where $\varphi$ is the Euler totient function and $\varphi(n)$ is the degree of $\Phi_n$). For $n\ge 3$, we use the notation: 
\[\Phi^{sym}_n(t)=t^{-\frac{\varphi(n)}{2}}\Phi_n(t)\]
For example $\Phi^{sym}_{10}(t)=t^{-2}-t^{-1}+1-t+t^2$
and $\Phi^{sym}_{14}(t)=t^{-3}-t^{-2}+t^{-1}-1+t-t^2+t^3$. Notice, that $\Phi^{sym}_n$ and $\Phi_n$ have the same roots.
One checks that the formula for $\Phit_{2m}(t)$ given above, $m$ odd,  simplifies to:
\[\Phit_{2m}(t):=\prod_{k|m,k>1}\Phi^{sym}_{2k}(t)\]

The paper is organized as follows: in section~\ref{sec:main} the main theorems are stated, while the notion of arrow diagrams and some proofs are postponed to two latter sections~\ref{sec:arrows} and \ref{sec:proofs}.

\section{Main results}\label{sec:main}

Let $W_{n,k}$, $n,k\in\Z$, $k\ge 0$, be the knot shown for $n=2$ and $k=3$ in Figure~\ref{fig:wnk}, in the form of an {\it arrow diagram}. In general, there are $n$ arrows on the left kink, and $k$ arrows arranged on $k$ strands, generalizing in an obvious way the case $k=3$ shown in this figure. When $n<0$, there are $|n|$ clockwise arrows on the left kink.
In short, the arrows correspond to fibers in the Hopf fibration of $S^3$. A detailed explanation of arrow diagrams is postponed to section~\ref{sec:arrows}. In section~\ref{sec:proofs}, we compute the Jones polynomials of the knots $W_{n,k}$, denoted $V_{W_{n,k}}$:

\begin{figure}[h]
\includegraphics{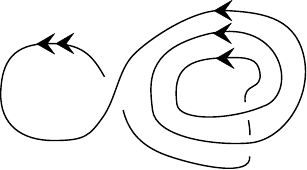}
\caption{$W_{2,3}$}
\label{fig:wnk}
\end{figure}

\begin{theorem}\label{thm:jonesWnk}
Let $n,k\in\Z$, $k\ge 0$. Then,
\begin{align*}
V_{W_{n,k}}=&\frac{t^{\frac{n(n-1)}{2}+k(k-1)-2nk}}{t^2-1}\left(-t^{(k+2)n+1}+t^{(k+1)(n+1)}(t^{k+1}+1)\right.\\
&\left.-t^{k(n+3)+1}+t-1\right)
\end{align*}
\end{theorem}

Denote by $D_{n,k}$ the terms in the big parenthesis in the formula for $V_{W_{n,k}}$ above. We want to check for which $n,k$ the polynomial $V_{W_{n,m}}$ is symmetric, i.e.  $V_{W_{n,m}}(t^{-1})=V_{W_{n,m}}(t)$.
It is easy to see that a necessary condion is that $D_{n,k}(t^{-1})=-t^aD_{n,k}(t)$ for some $a\in\Z$.
Say that such $D_{n,k}$ is {\it antipalindromic}.

For $k=0$, $W_{n,0}$ is an oval with $n\in\Z$ arrows on it. Such knots are torus knots, trivial if and only if $n\in\{1,0,-1,-2\}$, see~\cite{M1}. Thus, when $W_{n,0}$ is non trivial, its Jones polynomial is not symmetric.

\begin{theorem}\label{thm:cyclojones}
Suppose that $k>0$. The polynomial $V_{n,k}(t)$ is symmetric if and only if $n=k-1$, $k$, $2k$ or $2k+1$.
Furthermore, let $f(k)=k^2+k-1$ and $g(k)=2k^2-1$. Then, for $k>0$, 
\[V_{W_{k-1,k}}(t)=\Phit_{2f(k)}(t)\]
\[V_{W_{k,k}}(t)=\Phit_{2f(k+1)}(t)\]
\[V_{W_{2k,k}}(t)=V_{W_{2k+1,k}}(t)=\Phit_{2g(k+1)}(t)\]

\begin{proof}
We consider $D_{n,k}$. First we check when $t$ or $-1$ cancels with some other term.
The term $t$ cancels with a term $-t^a$, if $a=1$. This occurs when $n=0$ or $n=-3$.
The term $-1$ cancels if $n=-1$ or $n=-2$. One checks, that $D_{n,k}$ is not antipalindromic in all these cases,
except for $(n,k)=(0,1)$. That $W_{0,1}$ is trivial is very easy to check, see section~\ref{sec:arrows}. Notice, that
$\Phit_{2f(1)}=\Phit_2=1$.

Suppose now that $n$ and $k$ are such that neither $t$ nor $-1$ cancels, hence $n\ge 1$ or $n\le -4$.
Let $n_1=(k+2)n+1$, $n_2=k(n+3)+1$, $p_1=(k+1)(n+1)$ and $p_2=(k+1)(n+2)$ be the exponents of
the four remaining terms in $D_{n,k}$  (two negative and two positive ones).

Suppose that $n\le -4$. The four exponents are negative with the exception $n_2=0$ for $(n,k)=(-4,1)$. The highest terms are $t-1$ or $t-2$, the lowest is $-t^{n_1}$ and the gap between $n_1$ and any other exponent it at least $5$. Hence, $D_{n,k}$ is
not antipalindromic.

Suppose that $n\ge 1$. The four exponents are greater or equal to $4$. The four terms cannot all cancel out, since otherwise $V_{W_{n,k}}$ would not be a Laurent polynomial (one may also check case by case that, if a pair of terms cancels, another pair
does not cancel). Since there is a gap $1$ in the powers of $t$ and $-1$, in order for $D_{n,k}$ to be antipalindromic, it should contain another pair of terms $t^a-t^{a-1}$ for some $a\ge 5$. One has: 
\[p_1-n_1=k-n,\quad p_2-n_2=n-k+1,\quad p_2-n_1=2k-n+1,\quad p_1-n_2=n-2k\]
One of these four differences has to be equal to $1$, which gives four cases:
\begin{itemize}
\item $p_1-n_1=k-n=1$. Then, $p_2=n_2$ and these $2$ terms cancel out. Now:
\[D_{k-1,k}=t^{n_1}(t-1)+t-1=(t^{k^2+k-1}+1)(t-1)\]
Let $m=k^2+k-1$. One checks that:
\[V_{W_{k-1,k}}=\frac{t^{-\frac{m-1}{2}}}{t^2-1}(t^m+1)(t-1)=\Phit_{2m}=\Phit_{2f(k)}\]
\item $p_2-n_2=n-k+1=1$. Then, $p_1=n_1$ and:
\[D_{k,k}=t^{n_2}(t-1)+t-1=(t^{k^2+3k+1})(t-1)\]
Let $m=k^2+3k+1$. Again, one checks that:
\[V_{W_{k,k}}=\Phit_{2m}=\Phit_{2f(k+1)}\]
\item  $p_2-n_1=2k-n+1=1$. Then $p_1=n_2$ and:
\[D_{2k,k}=t^{n_1}(t-1)+t-1=(t^{2k^2+4k+1})(t-1)\]
Let $m=2k^2+4k+1$. One checks that:
\[V_{W_{2k,k}}=\Phit_{2m}=\Phit_{2g(k+1)}\]
\item $p_1-n_2=n-2k=1$. Then, $p_2=n_1$ and:
\[D_{2k+1,k}=t^{n_2}(t-1)+t-1=(t^{2k^2+4k+1})(t-1)\]
Let $m=2k^2+4k+1$. One checks that:
\[V_{W_{2k+1,k}}=\Phit({2m})=\Phit_{2g(k+1)}\]
\end{itemize}
\end{proof}
\end{theorem}

As an immedaite consequence, we get:

\begin{theorem}
There are infinitely many roots of unity that are zeros of Jones polynomials. Such roots are dense in the unit circle.
\begin{proof}
Since for any odd $m$, $\Phi_{2m}|\Phit_{2m}$, one has: for any $k>0$, $\zeta_{2f(k)}$, $\zeta_{2f(k+1)}$, $\zeta_{2g(k+1)}$ are zeros of some Jones polynomials.

Since $t^{2m}-1=(t^m-1)(t+1)\Phit_{2m}t^{\frac{m-1}{2}}$, the roots of $\Phit_{2m}$ are:
\[\zeta_{2m}, \zeta_{2m}^3\ldots\zeta_{2m}^{m-2},\zeta_{2m}^{m+2}\ldots \zeta_{2m}^{2m-1}\]
It is clear that the roots of the $\Phit_{2m}$'s, that are Jones polynomials, are dense in the unit circle, since there are infinitely many such $\Phit_{2m}$'s.
\end{proof}
\end{theorem}

The knots appearing in Theorem~\ref{thm:cyclojones} come in quadruplets for $k=1,2,3...$
In Table~\ref{tab:k4}, are shown the first four quadruplets, their Jones polynomials and the crossing numbers
(together with identification for knots up to $15$ crossings; also, the knot $W_{3,1}$ is $16n_{207543}$).

\begin{table}[h]
\begin{center}
\label{tab:k4}
\caption{$W_{n,k}$ with cyclotomic Jones polynomials up to $k=4$}
\begin{tabular}{c|c|c||c|c|c||c|c|c||c|c|c}
$K$&$V_K$&$c(K)$&$K$&$V_K$&$c(K)$&$K$&$V_K$&$c(K)$&$K$&$V_K$&$c(K)$\\
\hline
$W_{0,1}$ &  $1$& $0_1$ & $W_{1,1}$ &  $\Phi^{\scriptsize{sym}}_{10}$ & $4_1$ & 
$W_{2,1}$ & $\Phi^{sym}_{14}$ &  $9_{42}$ & $W_{3,1}$   & $\Phi^{sym}_{14}$ & $16$\\
\hline
$W_{1,2}$ & $\Phi^{sym}_{10}$ &$11n_{19}$ &  $W_{2,2}$ & $\Phi^{sym}_{22}$ & $\le 18$ & 
$W_{4,2}$ & $\Phi^{sym}_{34}$ & $\le 38$ & $W_{5,2}$ & $\Phi^{sym}_{34}$& $\le 52$\\
\hline
$W_{2,3}$ & $\Phi^{sym}_{22}$ &$\le 31$& $W_{3,3}$ & $\Phi^{sym}_{38}$  &$\le 43$& $W_{6,3}$ & $\Phi^{sym}_{62}$ &$\le 89$& $W_{7,3}$ & $\Phi^{sym}_{62}$ &$\le 108$\\
\hline &&&&&&&&&&&\\[-10pt]
 $W_{3,4}$ & $\Phi^{sym}_{38}$  &$\le 64$& $W_{4,4}$ & $\Phi^{sym}_{58}$  &$\le 79$& 
$W_{8,4}$ & $\Phit_{98}$ &$\le 159$& $W_{9,4}$ & $\Phit_{98}$ &$\le 184$
\end{tabular}
\end{center}
\end{table}

Notice that $98$ is the first index that is not twice a prime, hence $\Phit_{98}\neq\Phi^{sym}_{98}$.
The knots with $k>1$, except $W_{1,2}$, have more than $16$ crossings (since their Jones polynomials do not appear
up to $16$ crossings) and their crossing number seems to increase rapidely. From Lemma~\ref{lem:upperc} below, $c(W_{n,k})<=k^2+(n+k)^2-1$. Using Knotscape\cite{HT} (after removing the arrows in the diagrams, see section~\ref{sec:arrows}), the number of crossings can sometimes be reduced by $1$ or $2$. Knotscape handles diagrams up to $49$ crossings and allows
to check that the Alexander polynomial differentiates $W_{2,2}$ from $W_{2,3}$. Since $W_{5,2}$ has a diagram with $52$ crossings, its Alexander polynomial cannot be computed with Knotscape (in order to check whether it is different from $W_{4,2}$). Though it seems unlikely, it is possible that some $W_{2k,k}$ is the same knot as $W_{2k+1,k}$ and/or some $W_{k,k}$
is the same knot as $W_{k,k+1}$.

As an example, using Knotscape on a diagram of $W_{4,2}$ with $39$ crossings one gets a reduction to $38$ crossings
with the following DT code:

\texttt{38 1}\;\;\;\;\;\;\;\;\texttt{6 -14  16 -28  26 -42  68  40 -50 -60 -34  36 -44 -54  52  62 -20  46
-58  -4   2  56  64 -22  70 -76   8 -10 -66 -48  72 -74  24  12  38  18 -32  30 (best available reduction)}

We turn now to some properties of the knots $W_{n,k}$.

\begin{proposition}\label{prop:11}
The knots $W_{n,k}$ are $(1,1)$ knots. In particular they have tunnel number $1$, hence they
are prime.
\begin{proof}
We postpone the proof that $W_{n,k}$ are $(1,1)$ knots to section~\ref{sec:arrows}. Now $(1,1)$ knots have tunnel number $1$ (see~\cite{D}), hence they are prime (see~\cite{N, S}).
\end{proof}
\end{proposition}

\begin{proposition}
The knots $W_{n,k}$ with cyclotomic Jones polynomials are non alternating except for $0_1$ and $4_1$. There is no bound on the twist number of such knots.
\begin{proof}
Suppose that $W_{n,k}$ is non trivial and alternating with Jones polynomial equal to some $\Phit_{2m}$. From~\cite{DL}, its twist number equals $2$. Such knot can be either a connected sum of torus knots of type $(2,m)$ and $(2,n)$ or a 2-bridge knot. From Proposition~\ref{prop:11} $W_{k,n}$ is prime, so it has to be a 2-bridge knot. Using an explicit formula for Jones polynomials of 2-bridge knots with twist number $2$ in \cite{QYA}, one checks easily, that all such knots, except $4_1$, have Jones polynomials that are not equal to $\Phit_{2m}$ for any $m$. 

For the second part, it is shown in~\cite{CK2}, that for a family of links with cyclotomic Jones polynomials of unbounded span, there is no bound on the twist numbers of these links.
\end{proof}
\end{proposition}

We turn now to some obstructions for roots of unity being zeros of Jones polynomials.  

It is well known that the Jones polynomial has special values in $1$, $\zeta_3$, $i$ and $\zeta_6$, see~\cite{J, LM}.
For a knot $K$, $V_{K}(1)=V_{K}(\zeta_3)=1$, $V_K(\zeta_4)=\pm 1$ and $V_K(\zeta_6)=\pm(i\sqrt{3})^n$, $n\in\No$.
This allows to exclude some roots of unity as zeros of Jones polynomials:

\begin{theorem}\label{thm:excluded}
For $k\in\No$, let $N=p^k$, $3p^k$, $4p^k$, with $p$ prime; or $N=6p^k$  with $p\neq 3$ prime. 
Then $\Phi_N$ cannot divide any Jones polynomial.
\begin{proof}
We check the values of cyclotomic polynomials in $1$, $\zeta_3$, $i$ and $\zeta_6$. 
From \cite{BHM}, we have for $p$ prime:
\[\Phi_{p^k}(1)=p$ for $k>0$ (and $\Phi_1(1)=0)\]
\[|\Phi_{3p^k}(\zeta_3)|=\Phi_{4p^k}(i)=|\Phi_{6p^k}(\zeta_6)|=p\]

We see that none of these polynomials can divide a Jones polynomial $V_K$, since $V_K(1)=V_K(\zeta_3)=|V_K(i)|=1$ and
$|V_K(\zeta_6)|=\sqrt{3}^n$, except if $p=3$ in the case of $\Phi_{6p^k}$, which we excluded in our assumptions.
\end{proof}
\end{theorem}

We can also exclude easily some $\Phit_{2k}$ as divisors of Jones polynomial:

\begin{proposition}
Let $k$ be odd. If $\Phit_{2k}$ divides a Jones polynomial, then $k$ is not divisible by $3$.
\begin{proof}
If $3|k$, then $\Phi_6|\Phit_{2k}$, hence $\Phit_{2k}(\zeta_6)=0$, which is impossible for a divisor of a Jones polynomial. 
\end{proof}
\end{proposition}

Notice, that all roots of unity appearing as zeros of Jones polynomials in Theorem~\ref{thm:cyclojones} are of the form $\zeta_{2k}$, with $k$ odd, $3\nmid k$. It is natural to ask what other roots of unity $\zeta_N$ can be zeros of Jones polynomials of knots.
Using Theorem~\ref{thm:excluded}, the smallest possible $N$ for such roots are $18,26,35,40,45,46,50,54,55,56,60$. Let us sum this up:

\begin{question}
Is there a knot with Jones polynomial having a zero in $\zeta_N$ such that:
$4|N$; $N$ is odd; $3|N$; or $N=2k$, $k$ odd, $3\nmid k$ but $N$ not coming from Theorem~\ref{thm:cyclojones}?
\end{question}

One may also ask the question, whether there are infinitely many primes such that $\Phi^{sym}_{2p}$ are Jones polynomials of some knots. A positive answer would follow, if there were infinitely many primes in the image of $f$ or $g$ from Theorem~\ref{thm:cyclojones} (two special cases of the Bunyakovsky conjecture). 

Since a Mersenne prime $2^p-1$, with $p>2$, satisfies $2^p-1=2(2^\frac{p-1}{2})^2-1=g(2^\frac{p-1}{2})$, we get:

\begin{corollary}
Let $N=2^p-1$, $p>2$, be a Mersenne prime. Then $\Phi^{sym}_{2N}$ is the Jones polynomial of a knot.
\end{corollary}

\section{Arrow diagrams}\label{sec:arrows}

\begin{figure}[h]
\centering
\includegraphics{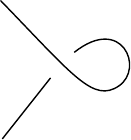} \raisebox{2 em}{$\; \underset{\longleftrightarrow}{\Omega_1} \;$} \includegraphics{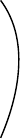}
\hspace{2em}
\includegraphics{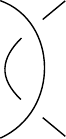} \raisebox{2 em}{$\; \underset{\longleftrightarrow}{\Omega_2} \;$} \includegraphics{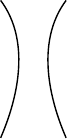}
\hspace{2em}
\includegraphics{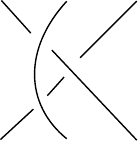} \raisebox{2 em}{$\; \underset{\longleftrightarrow}{\Omega_3} \;$} \includegraphics{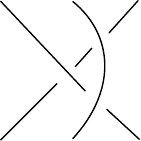}
\hspace{2em} \\[2em]
\includegraphics{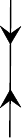} \raisebox{2 em}{$\; \underset{\longleftrightarrow}{\Omega_4} \;$} \includegraphics{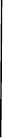}
\raisebox{2 em}{$\; \underset{\longleftrightarrow}{\Omega_4} \;$} \includegraphics{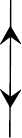} \hspace{5em}
\includegraphics{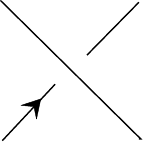} \raisebox{2 em}{$\; \underset{\longleftrightarrow}{\Omega_5} \;$} \includegraphics{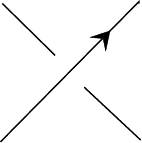}
\\[2em]
\includegraphics{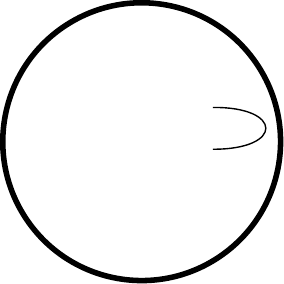} \raisebox{4.2 em}{$ \;
 \underset{\longleftrightarrow}{\Omega_{\infty}} \; $ } \includegraphics{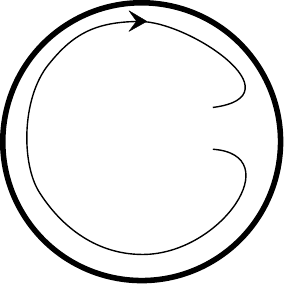}
\caption{Reidemeister moves}\label{fig:reid}
\end{figure}

Arrow diagrams where introduced in~\cite{MD} for links in $F\times S^1$, where $F$ is an orientable surface. They were subsequently extended for links in Seifert manifolds (see~\cite{GM, MM1,MM2}). In~\cite{M1}, they were applied for links in $S^3$: it was shown there, that projections of links under the Hopf fibration from $S^3$ to $S^2$ can be encoded with arrow diagrams in a disk: such a diagram is like a usual diagram of a link, except that it is in a disk and there may be some arrows on it, outside crossings. Two arrow diagrams represent the same link if and only if one diagram can be transformed into the other with a series of six Reidemeister moves, see Figure~\ref{fig:reid}. For the $\Om_\infty$ move in this figure, the boundary of the disk is drawn in thick. For simplicity we can also omit this boundary when picturing arrow diagrams (as we have done in Figure~\ref{fig:wnk}). 

A detailed interpretation of the arrow diagrams and Reidemeister moves can be found in~\cite{M1}. One can picture easily a link $L$ from its arrow diagram $D$ in the following way: pick a solid torus $T=C\times S^1$, $C$ a disk, consisting of some oriented fibers $p\times\ S^1$, $p\in C$, in the Hopf fibration of $S^3$. Let $S^1=I\cup I'$ consist of two intervals glued along their endpoints.
Then $T=B\cup B'$, where $B=C\times I$ and $B'=C\times I'$ are two balls. If there are no arrows in $D$, $L$ lies entirely in $B$. Otherwise it lies in $B$ except for some neighborhoods of the arrows where it goes through $B'$ along an oriented fiber and the orientation of the arrow agrees with the orientation of the fiber.

We turn now to the proof of Proposition~\ref{prop:11}. We want to show that the knots $W_{n,k}$ are $(1,1)$ knots. Recall from \cite{D} that a link $L$ admits a $(g,b)$ decomposition, if there is a genus $g$ Heegard splitting $(V_0,V_1)$ of $S^3$ such that $V_i$ intersects $L$ in $b$ trivial arcs, for $i\in\{0,1\}$. To show that $W_{n,k}$ is a $(1,1)$ knot, we need to show that it intersects each $T_i$ in a trivial arc, for a Heegard spliting of $S^3$ into two solid tori $T_i$, $i\in\{0,1\}$.

\begin{figure}[h]
\includegraphics{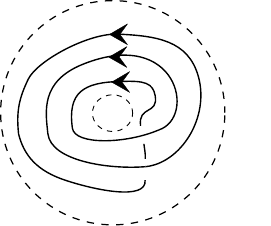}
\caption{An annulus monotonic diagram, drawn inside the annulus}
\label{fig:annulus}
\end{figure}

We say that an arrow diagram of a knot is {\it annulus monotonic}, if there is an annulus $A=S^1\times I$, containing the diagram and such that the curve of the diagram has exactly one minimum and one maximum w.r.t. $I$. 
Applying $\Om_\infty$ on the left kink of $W_{n,k}$ (see Figure~\ref{fig:wnk}), we obtain a diagram consiting of a spiral with some arrows on it. Such a diagram is clearly annulus monotonic, see Figure~\ref{fig:annulus}. Proposition~\ref{prop:11} now follows directly from the following:

\begin{lemma}
Suppose that a knot $K$ has an annulus monotonic arrow diagram $D$. Then $K$ is a $(1,1)$-knot.
\begin{proof}
Let $A=S^1\times I$ be an annulus containing $D$ and such that $D$ has exactly one minimum and one maximum w.r.t. $I$. 
The closure of $S^3\setminus (A\times S^1)$ consists of two solid tori $T_0$ and $T_1$, chosen so that $T_i\cap (A\times S^1)=(S^1\times {i})\times S^1$, $i\in\{0,1\}$.

Cut $I$ into $I_0=[0,a]$ and $I_1=[a,1]$ for some $a\in (0,1)$, so that the intersection $(S^1\times I_0)\cap D$ is a small trivial arc. $A$ decomposes into two annuli $A_i=S^1\times I_i$, $i\in\{0,1\}$. Let $T'_i=T_i\cup (A_i\times S^1)$ be two solid tori, $i\in\{0,1\}$, so that $S^3=T'_0\cup T'_1$. Then $T'_0\cap K$ is clearly a trivial arc in $T'_0$. We claim that $T'_1\cap K$ is also a trivial arc in $T'_1$. Let $I_2=[b,1]$, $b>a$, be such that $(S^1\times I_2)\cap D$ is a small trivial arc. Let $A_2=S^1\times I_2$. Since $D$ is annulus monotonic, the pair $(A_1\times S^1,K\cap (A_1\times S^1))$ can be isotoped
to $(A_2\times S^1,K\cap (A_2\times S^1))$, by removing the tori $(S^1\times {c})\times S^1$ for $c$ from $a$ to $b$. Such isotopy clearly extends to $(T'_1,K\cap T'_1)$, so $K\cap T'_1$ is a trivial arc in $T'_1$. Thus $K$ is a $(1,1)$ knot.
\end{proof}
\end{lemma}

We remark here, that for some $n$'s, hypothetical knots with Jones polynomials $\Phi^{sym}_n$ would only admit a $(g,b)$ decomposition with large $g+b$.
Indeed, in \cite{BHM} the values of cyclotomic polynomials in $\zeta_5$ are computed. It is shown there, 
that $|\Phi_n(\zeta_5)|$ can be arbitrarily large for some $n$'s. For example it grows very fast with the number of primes in the decomposition of $n$, when $n$ is a product of an odd number of distinct primes congruent to $2$ or $3$ modulo $5$. For instance, for $n=2\;3\;7\;13\;17$, one checks that this module is approximately $2207$.
On the other hand, it follows from \cite{MSY} that, if a knot $K$ admits a $(g,b)$ decomposition, its Jones polynomial $V_K$ satisfies $|V_K(\zeta_5)|\le \alpha^g\beta^{b-1}$, where $\alpha>1$ and $\beta>1$ can be explicitely computed. Hence, a large module implies large $g+b$.

It was shown in~\cite{M1}, that the usual blackboard framing for links obtained from their diagrams extends to arrow diagrams and that such framing is invariant under all Reidemeister moves except $\Om_1$. In particular, to compute the writhe of a framed link represented by an arrow diagram, one may eliminate all arrows without using $\Om_1$, then sum the signs of all crossings in the arrowless diagram.

We present now a formula for the writhe of any arrow diagram of a knot. This formula holds also for oriented links.

Let $D$ be an oriented arrow diagram. Let $r$ be an arrow in $D$.  The {\it sign} of $r$, denoted $\epsilon(r)$, is defined as follows: $\epsilon(r)=1$ (resp. $\epsilon(r)=-1$), if $r$ points in the same (resp. opposite) direction as the orientation of the diagram. We also say that $r$ is {\it positive} (resp. {\it negative}).
The winding number of $r$, denoted $ind(r)$, is by definition the winding number $ind_D(P)$, where $D$ is the diagram considered as an oriented curve and $P$ is a point close to $r$, to the right of $D$ according to the orientation of $D$.

For example, consider $W_{2,3}$ in Figure~\ref{fig:wnk}. Orient it so that the left kink is oriented clockwise. Then the $3$
arrows on the right are positive and the two arrows on the left are negative. Also the winding numbers of the arrows on the right are $0$, $1$ and $2$, wheras the $2$ arrows on the left have winding number $-1$.

Denote by $w(D)$ the writhe of the framed knot represented by the arrow diagram $D$. Denote by $\bar{w}(D)$ the writhe, when all arrows in $D$ are ignored (it is sum of the signs of crossings in $D$). We have the following formula for the writhe:

\begin{lemma}\label{lem:wr_formula}
Let $D$ be an oriented arrow diagram. Let $n=\Sum_{r}\epsilon(r)$, the sum taken over all arrows of $D$.
Then:

\[w(D)=\bar{w}(D)+\Sum_r 2\epsilon(r)ind(r)+n(n+1)\]
\begin{proof}
We remove with Reidemeister moves all arrows in $D$ keeping track of the signs of the crossings that appear. We do not use $\Om_1$, thus the writhe is unchanged.

Consider an arrow $r$ in $D$. We push it next to the boundary of the diagram in such a way that the orientation of the arc next to the arrow agrees with the counterclockwise orientation of the boundary of the diagram (see Figure~\ref{fig:pn_arrows} (left),
where $3$ arrows have been pushed and the arcs are oriented as wished). 

\begin{figure}[h]
\includegraphics{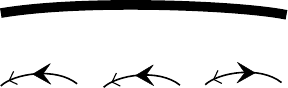}\qquad\qquad\includegraphics{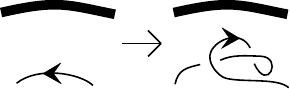}
\caption{Two positive, one negative arrow (left); pushing an arrow through the arc next to it (right)}
\label{fig:pn_arrows}
\end{figure}

To achieve this, we use $\Om_2$ and $\Om_5$ moves repeatedly. When $r$ crosses an arc, two positive or two negative crossings appear. One checks that the total contribution, when $r$ is next to the boundary, is $2\epsilon(r)ind(r)$. Notice that when $r$ is next to the boundary, but the orientation of the arc is not the desired one, then $ind(r)=-1$ and $r$ has to be pushed once through a piece of arc next to it, see Figure~\ref{fig:pn_arrows} (right).
After all the arrows have been pushed, so they are as in Figure~\ref{fig:pn_arrows} (left), the sum of the signs of all crossings is $\bar{w}(D)+\Sum_r 2\epsilon(r)ind(r)$.

Suppose now, that there are $a$ positive arrows and $b$ negative ones, so that $n=a-b$. Push every positive arrow through an arc next to it as in Figure~\ref{fig:pn_arrows} (right). This adds $2a$ positive crossings. Now any arrow $r$ can be eliminated with $\Om_\infty$ followed by $\Om_4$. We push the remaining arrows through the arc created by $\Om_\infty$. One checks that if an arrow $r'$ is pushed through the arc coming from $r$, this adds two positive (resp. negative) crossings if $\epsilon(r)=\epsilon(r')$ (resp. $\epsilon(r)=-\epsilon(r')$). Then one repeats the process with the arrow $r'$ (eliminating it and pushing all other arrows through it). Hence, at the end any pair of arrows $r$ and $r'$ contributes $2\epsilon(r)\epsilon(r')$ to the writhe. The total contribution to the writhe of this second part is thus:
\[2a+a(a-1)+b(b-1)-2ab=(a-b)(a-b+1)=n(n+1)\]
Combined with the first part, this gives the required formula.

\end{proof}
\end{lemma}

Applying Lemma~\ref{lem:wr_formula} to $W_{n,k}$ we get:

\begin{lemma}\label{lem:writheWL}
Let $W_{n,k}$ stand for the diagram in Figure\ref{fig:wnk}, as well as for the framed
knot represented by this diagram. Then:
\[w(W_{n,k})=n^2+n+2k^2+k-2nk\]

\begin{proof}
Orient $W_{n,k}$ so that the $k$ arrows are positive. If $n>0$ then the $n$ arrows are negative. If $n<0$ then the $|n|$ arrows are positive. The $k$ arrows have winding numbers $0, 1, 2,\ldots,k-1$. The
$n$ arrows have all winding number $-1$. Also, $\bar{w}(W_{n,k})=k$. Hence:
\[w(W_{n,k})=k+2(-n)(-1)+2(1+2+\ldots+k-1)+(k-n)(k-n+1)\]
\[=k+2n+k(k-1)+(k-n)(k-n+1)=n^2+n+2k^2+k-2nk\]
\end{proof}
\end{lemma}

For an upper estimate of the number of crossings, $c(W_{n,k})$, we use Lemma~1 from~\cite{M1}. It states that, if a diagram $D$ has $k$ crossings and all its arrows are next to the boundary, with $a$ of them removable (i.e. one can remove them with $\Om_\infty$ followed by $\Om_4$) and $b>0$ of them non removable, then $c(K)\le k+b-1+(a+b)(a+b-1)$. We get:

\begin{lemma}\label{lem:upperc}
For $n\ge 0$, $k\ge 0$ and $k+n>0$, one has:
\[c(W_{n,k})\le k^2+(n+k)^2-1\]
\begin{proof}
Starting with the diagram of $W_{n,k}$ shown in Figure~\ref{fig:wnk}, we push $k-1$ arrows so that they are next to the boundary. We get a diagram with $k+2(1+2+\ldots+k-1)=k+k(k-1)=k^2$ crossings. Then, we can apply Lemma~1 from~\cite{M1}. Since $n\ge 0$, all arrows will be non removable. Since $n+k>0$,  there is at least one non removable arrow. Thus, $c(W_{n,k})\le k^2+(n+k)-1+(n+k)(n+k-1)=k^2+(n+k)^2-1$.
\end{proof}
\end{lemma}

We end this section with a visualization of any knot $W_{n,k}$. Such knot is obtained by a small modification from a pair of torus knots lying on the boundary of a thickened Hopf link. It was shown in~\cite{M1} how to get some simple arrow diagrams of torus knots: one checks that $W_{n,0}$ is the torus knot $T(n,n+1)$ and $W_{0,k}$ is the torus knot $T(k,2k+1)$. Consider the diagram of $W_{n,k}$ in Figure~\ref{fig:wnk}. Let $W^s_{n,k}$ be the diagram of a $2$-component link, obtained from $W_{n,k}$ by smoothing vertically the crossing next to the $n$ arrows. The components of $W^s_{n,k}$ are torus knots $T(n,n+1)$ and $T(k,2k+1)$.
Let $D$ and $D'$ be two disjoint disks, such that $D$ contains the $n$ arrows, $D'$ contains the $k$ arrows and $W^s_{n,k}$ is contained in $D\cup D'$. Let $T$, resp. $T'$, be two solid tori consisting of fibers intersecting $D$, resp. $D'$, in the Hopf fibration of $S^3$. Then $T$ and $T'$ form a thickened Hopf link. The torus $(n,n+1)$ component of $W^s_{n,k}$ can be pushed onto $\partial T$ and the torus $(k,2k+1)$ component can be pushed onto $\partial T'$. Then $W_{n,k}$ is obtained from such two linked torus knots by reverting the smoothing back to the crossing.

\section{Jones polynomials of the knots $W_{n,k}$}\label{sec:proofs}

\begin{figure}[h]
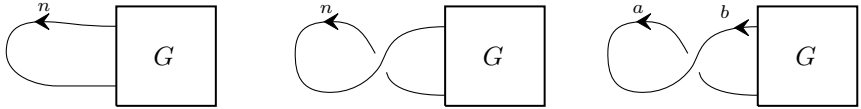

\oversimage{Gn}{14,56}{\scriptsize $n$}{65,32}{$G$}
\qquad\oversimage{Gpn}{10,50}{\scriptsize $n$}{75,29}{$G$}
\qquad\overssimage{Gpab}{10,50}{\scriptsize $a$}{45,48}{\scriptsize$b$}{75,29}{$G$}
\caption{$G_n$, $G'_n$ and $G'_{a,b}$}
\label{fig:gn}
\end{figure}

Let $G_n$, $G'_n$ and $G'_{a,b}$, $n,a,b\in\Z$, be the arrow diagrams shown in Figure~\ref{fig:gn}. In the box is an arrow tangle $G$ (a tangle with, possibly, some arrows on it). Let $g_n$, $g'_n$ and $g'_{a,b}$ be the Kauffman brackets of, respectively, $G_n$, $G'_n$ and $G'_{a,b}$.
We want to express $g'_n$ with some $g_k$'s. In order to do it, we will use $g'_{a,b}$'s.

It is useful to define for $n\ge0$ the sum:
\[S_n=A^{n}g_{-n}+A^{n-2}g_{-n+2}+\ldots+A^{-n+2}g_{n-2}+A^{-n}g_n=\Sum_{i=0}^n A^{n-2i}g_{-n+2i}\]
Extend $S_n$ for negative $n$, by defining $S_{-1}=0$ and, for $n<-1$:
\[S_n=-S_{|n|-2}\]

\begin{lemma}\label{lem:kink}
For $n\in\Z$:
\[g'_n=(A^{-1}-A^3)A^nS_n-A^{2n-1}g_{-n}\]

\begin{proof}
One checks easily that the formula holds for $n=0$ and $n=-1$.

From the defining relations of the Kauffman bracket, we get:

\[
\raisebox{-7pt}{\includegraphics{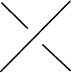}}=A^2 \raisebox{-7pt}{\includegraphics{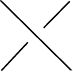}}
+(A^{-1}-A^{3}) \raisebox{-7pt}{\includegraphics{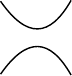}}\]

Using this relation and $\Om_5$ and $\Om_4$ moves, we get:
\[g'_{a,b}=A^2g'_{a-1,b-1}+(A^{-1}-A^3)g_{a+b}\label{eq:ttt}\tag{*}\]

Suppose that $n\ge 1$. Iterating equation~(\ref{eq:ttt}) until $g'_{0,-n}=-A^3g_{-n}$, we get:
\begin{align*}
g'_{n,0}=&\quad A^2g'_{n-1,-1}+(A^{-1}-A^3)g_{n}\\
=&\quad A^4g'_{n-2,-2}+A^2(A^{-1}-A^3)g_{n-2}+(A^{-1}-A^3)g_{n}\\
=&\quad \ldots\\
= &\quad -A^3A^{2n}g_{-n}+(A^{-1}-A^3)\left(g_n+A^2g_{n-2}+\ldots+A^{2n-2}g_{-n+2}\right)\\
=&\quad -A^{2n-1}g_{-n}+(A^{-1}-A^3)\left(g_n+A^2g_{n-2}+\ldots+A^{2n}g_{-n}\right)\\
=&\quad (A^{-1}-A^3) A^n S_n -A^{2n-1}g_{-n}
\end{align*}

Suppose now that $n\le -2$. Rewriting equation~(\ref{eq:ttt}) and replacting $a$ by $a+1$ and $b$ by $b+1$ one gets:
\[g'_{a,b}=A^{-2}g'_{a+1,b+1}-A^{-2}(A^{-1}-A^3)g_{a+b+2}\]
Iterating until $g'_{n+|n|,|n|}=g'_{0,-n}=-A^3g_{-n}$, we get:
\begin{align*}
g'_{n,0}=&\quad A^{-2}g'_{n+1,1}-A^{-2}(A^{-1}-A^3)g_{n+2}\\
=&\quad A^{-4}g'_{n+2,2}-A^{-4}(A^{-1}-A^3)g_{n+4}-A^{-2}(A^{-1}-A^3)g_{n+2}\\
=&\quad \ldots\\
= &\quad -A^3A^{2n}g_{-n}-A^{-2}(A^{-1}-A^3)\left(g_{n+2}+A^{-2}g_{n+4}+\ldots+A^{2n+2}g_{-n}\right)\\
= &\quad -A^{2n-1}g_{-n}-A^{-2}(A^{-1}-A^3)\left(g_{n+2}+A^{-2}g_{n+4}+\ldots+A^{2n+4}g_{-n-2}\right)\\
= &\quad -A^{2n-1}g_{-n}-A^{-2}(A^{-1}-A^3)A^{n+2}\left(A^{-n-2}g_{n+2}+\ldots+A^{n+2}g_{-n-2}\right)\\
= &\quad -A^{2n-1}g_{-n}-(A^{-1}-A^3)A^{n}S_{|n|-2}\\
= &\quad  (A^{-1}-A^3)A^{n}S_{n}-A^{2n-1}g_{-n}
\end{align*}

\end{proof}
\end{lemma}

We now prove Theorem~\ref{thm:jonesWnk} by induction on $k$. $W_{n,0}$ is an oval with $n\in\Z$ arrows on it. This is the torus knot $T(n,n+1)$ if $n\ge 0$ and, for $n<0$, $W_{n,0}=W_{-1-n,0}$ (use $\Om_\infty$ and $\Om_4$ moves). One checks that for $k=0$ the formula in Theorem~\ref{thm:jonesWnk} is the correct formula for such torus knots (see also~\cite{M1}).

We restate Theorem~\ref{thm:jonesWnk} in terms of the Kauffman bracket using the formula $V_K(t)=(-A)^{-3w(K)}<K>$,
where $w(K)$ is the writhe of $K$ and $t=A^{-4}$.
From Lemma~\ref{lem:writheWL}, we have $w(W_{n,k})=n^2+n+2k^2+k-2nk$, hence $(-1)^{3w(W_{n,k})}=(-1)^k$ and:
\[<W_{n,k}>=(-1)^kA^{3(n+1)n+ 3k(2k+1-2n)}V_{W_{n,k}}\]

One checks, that Theorem~\ref{thm:jonesWnk} can be restated as:

\begin{proposition}
\[<W_{n,k}>(A^{-8}-1)=(-1)^kA^{n^2-2kn+2k^2+n-k-8}\]
\[\left((1-A^4)A^{4kn+4n+8k+4}-A^{4n-4k+4}+A^{4k+4}-A^{8k-4n+4}+1\right)\]

\begin{proof}
For $k=0$ the formula is correct, since it is correct for the Jones polynomial and we just restate it with the Kauffman bracket using the writhe.

Let $k\ge 0$. Assume that the formula holds for $<W_{n,k}>$, $n\in\Z$.
We use Lemma~\ref{lem:kink}, with $W_{n,k+1}=G'_n$. One has to identify $G_a$ (in that lemma), for any $a\in\Z$: it has a diagram shown in Figure~\ref{fig:WGn} (for $k=2$ as an example). Using a single $\Om_\infty$ move on the strand with the $a+1$ arrows, one gets $W_{-a-2,k}$ (one extra arrow comes from the move). Since this move does not change the writhe, $<G_a>=<W_{-a-2,k}>$.

\begin{figure}[h]
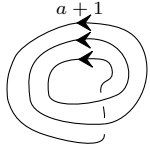

\overimage{WGn}{35,92}{\scriptsize $a+1$}
\caption{Identifying $G_a$}
\label{fig:WGn}
\end{figure}

Recall the notations used in Lemma~\ref{lem:kink}:
\[g'_n=<G'_n>,\quad g_a=<G_a>,\quad S_n=\Sum_{i=0}^n A^{n-2i}g_{-n+2i}\]
Also, by definition: $S_{-1}=0$ and $S_n=-S_{|n|-2}$ for $n<-1$.

Let:
\begin{align*}
S'_{n}=&(-1)^kA^{n^2-2kn-6n+2k^2-k-10}\left(-A^{4kn+12n+12k+16}\right.\\
&\left.+A^{4kn+8n+4k}+A^{4n+8k+8}-A^{8n}\right)
\end{align*}
One checks that $S'_{-1}=0$ and $S'_n=-S'_{-n-2}$ for any $n\in\Z$.

We claim that for any $n\in\Z$:
\[S_{n}(A^{-8}-1)=S'_{n} \tag{**}\label{ssp}\]

Because of the skew-symmetry of both $S_n$ and $S'_n$ around $-1$, it is sufficient to prove~(\ref{ssp}) for $n\ge -1$.

It is true for $n=-1$. Now $S_{0}=g_0=<W_{-2,k}>$. By induction on $k$:
\[<W_{-2,k}>(A^{-8}-1)=(-1)^k A^{2k^2-k-10}(-A^{12k+16}+A^{8k+8}+A^{4k}-1)=S'_{0}\]

One has obviously:
\begin{align*}
S_{n+2}&=S_n+A^{n+2}g_{-n-2}+A^{-n-2}g_{n+2}\\
&=S_n+A^{n+2}<W_{n,k}>+A^{-n-2}<W_{-n-4,k}>
\end{align*}

Thus, to prove~(\ref{ssp}), we need to show
that:
\[S'_{n+2}=S'_n+(A^{-8}-1)\left(A^{n+2}<W_{n,k}>+A^{-n-2}<W_{-n-4,k}>\right)\]
One checks:
\[S'_{n+2}=(-1)^kA^{n^2-2kn+2k^2-2n-5k-18}\left(-A^{4kn+12n+20k+40}\right.\]
\[\left.+A^{4kn+8n+12k+16}-A^{8n+16}+A^{4n+8k+16}\right)\]
By induction on $k$:
\begin{align*}
&S'_n+(A^{-8}-1)\left(A^{n+2}<W_{n,k}>+A^{-n-2}<W_{-n-4,k}>\right)=\\
&(-1)^kA^{n^2-2kn-6n+2k^2-k-10}\left(-A^{4kn+12n+12k+16}\right.\\
&\left.+A^{4kn+8n+4k}+A^{4n+8k+8}-A^{8n}\right)+(-1)^kA^{n^2-2kn+2k^2+2n-k-6}\\
&\left((1-A^4)A^{4kn+4n+8k+4}-A^{4n-4k+4}+A^{4k+4}-A^{8k-4n+4}+1\right)\\
&+(-1)^k A^{n^2+2kn+2k^2+6n+7k+2}\left((1-A^4)A^{-4kn-4n-8k-12}-A^{-4n-4k-12}\right.\\
&\left.+A^{4k+4}-A^{8k+4n+20}+1\right)=(-1)^kA^{n^2-2kn+2k^2-2n-5k-18}\\
&\left(-A^{4kn+8n+16k+24}+A^{4kn+4n+8k+8}+A^{12k+16}-A^{4n+4k+8}\right.\\
&+(1-A^4)A^{4kn+8n+12k+16}-A^{8n+16}+A^{4n+8k+16}-A^{12k+16}+A^{4n+4k+12}\\
&+(1-A^4)A^{4n+4k+8}-A^{4kn+4n+8k+8}+A^{4kn+8n+16k+24}\\
&\left.-A^{4kn+12n+20k+40}+A^{4kn+8n+12k+20}\right)=S'_{n+2}
\end{align*}

Since $g_{-n}=<W_{n-2,k}>$, from Lemma~\ref{lem:kink} we get:
\[<W_{n,k+1}>=(A^{-1}-A^3)A^{n}S_{n}-A^{2n-1}<W_{n-2,k}>\]
Hence:
\begin{align*}
&<W_{n,k+1}>(A^{-8}-1)=(A^{-1}-A^3)A^{n}S'_{n}-A^{2n-1}<W_{n-2,k}>(A^{-8}-1)\\
&=(-1)^k(A^{-1}-A^3)A^{n^2-2kn-5n+2k^2-k-10}\left(-A^{4kn+12n+12k+16}+A^{4kn+8n+4k}\right.\\
&\left.+A^{4n+8k+8}-A^{8n}\right)-(-1)^kA^{n^2-2kn+2k^2-n+3k-7}\left((1-A^4)A^{4kn+4n-4}\right.\\
&\left.-A^{4n-4k-4}+A^{4k+4}-A^{8k-4n+12}+1\right)=(-1)^kA^{n^2-2kn+2k^2-n+3k-7}\\
&\left(-(1-A^4)A^{4kn+8n+8k+12}+(1-A^4)A^{4kn+4n-4}+(1-A^4)A^{4k+4}\right.\\
&-(1-A^4)A^{4n-4k-4}-(1-A^4)A^{4kn+4n-4}+A^{4n-4k-4}-A^{4k+4}\\
&\left.+A^{8k-4n+12}-1\right)=(-1)^{k+1}A^{n^2-2kn+2k^2-n+3k-7}\\
&\left((1-A^4)A^{4kn+8n+8k+12}-A^{4n-4k}+A^{4k+8}-A^{8k-4n+12}+1\right)\\
&=(-1)^{k+1}A^{n^2-2(k+1)n+2(k+1)^2+n-(k+1)-8}\left((1-A^4)A^{4(k+1)n+4n+8(k+1)+4}\right.\\
&\left.-A^{4n-4(k+1)+4}+A^{4(k+1)+4}-A^{8(k+1)-4n+4}+1\right)
\end{align*}


Thus the formula holds for $<W_{n,k+1}>$ and we are done.
\end{proof}
\end{proposition}

\end{document}